\documentclass[MSNbibl,number,citesort,seceqn,dvips]{arxbj}
\usepackage{upgreek}
\usepackage{graphicx}

\aid{0}
\volume{20}
\issue{4}
\pubyear{2014}
\firstpage{1979}
\lastpage{1998}
\doi{10.3150/13-BEJ547} 
\setattribute{copyright}{owner}{In the Public Domain}

\makeatletter
\newtheorem{prp}{Proposition}[section]
\newtheorem{them}{Theorem}[section]
\newtheorem{lem}{Lemma}[section]

\renewcommand{\epsilon}{\varepsilon}
\newcommand{\xt}{\hat{\mu}_{\mathrm{opt}}}
\newcommand{\tx}{\hat{\mu}_{\mathrm{plug}}}
\newcommand{\w}{\omega}
\newcommand{\h}{\rho}
\renewcommand{\citep}{\cite}
\makeatother

\begin{document}
\begin{frontmatter}

\title{Restricted likelihood representation and decision-theoretic
aspects of meta-analysis}
\runtitle{Decision-theoretic aspects of meta-analysis}

\begin{aug}
\author{\inits{A.L.}\fnms{Andrew L.} \snm{Rukhin}\corref{}\ead[label=e1]{andrew.rukhin@nist.gov}}
\address{National Institute of Standards and Technology,
100 Bureau Dr., Gaithersburg, MD 20899, USA.\\ \printead{e1}}
\end{aug}

\received{\smonth{1} \syear{2013}}
\revised{\smonth{4} \syear{2013}}

%
\begin{abstract}
In the random-effects model of meta-analysis a canonical representation
of the restricted likelihood function is
obtained. This representation relates the mean effect and the
heterogeneity variance estimation problems.
An explicit form of the variance of weighted means statistics
determined by means of a quadratic form
is found. The behavior of the mean squared error for large
heterogeneity variance is elucidated.
It is noted that the sample mean is not admissible nor minimax under a
natural risk function for the number of
studies exceeding three.
\end{abstract}

%

\end{frontmatter}

\section{Parameter estimation in meta-analysis: Random-effects model}
\label{in}

In the simplest random-effects model of meta-analysis involving, say,
$n$ studies the data is supposed to consist
of treatment effect estimators $x_{i}, i=1, \ldots, n$, which have the form
\[
x_{i}= \mu+b_i+ \epsilon_i.
\]
Here $\mu$ is an unknown common mean, $b_i$ is zero mean between-study
effect with variance
$\tau^2$, $\tau^2 \geq0$, and $\epsilon_i$ represents the
measurement error of the $i$th study,
with variance $\sigma_i^2, \sigma_i^2 >0$. Then the variance of $x_i$
is $\tau^2+\sigma^2_i$.
In practice $\sigma_i$ is often treated as a given constant,
$s_i$, which is the reported standard error or uncertainty of the
$i$th study.

The considered here problem is
that of estimation of the common mean $\mu$ and of the heterogeneity
variance $\tau^2$ from the statistical
decision theory point of view under normality assumption. If
$\tau^2$ is known, then the best unbiased estimator of $\mu$ is the
weighted means
statistic, $\xt=\sum\w_i^0 x_i$, with the normalized weights,
%
\begin{equation}
\label{we} \w_i^0 = \frac{1}{\tau^2+s_i^2} \biggl(\sum
_k \frac{1}{ \tau
^2+s_k^2 } \biggr) ^{-1},
\end{equation}
$\sum\w_i^0=1$. Its variance has the form
\[
\operatorname{Var}(\xt)= \biggl[\sum_i
\frac{1}{ \tau^2+s_i^2 } \biggr] ^{-1}.
\]
When $\tau^2$ is unknown, to estimate $\mu$ the common practice uses
a plug-in version of $\xt$,
%
\begin{equation}
\label{es} \tx=\sum_i\frac{x_i}{\hat{\tau}^2+s_i^2}
\biggl(\sum_i\frac
{1}{\hat{\tau}^2+s_i^2 } \biggr)
^{-1},
\end{equation}
so that an estimator $\hat{\tau}^2$ of $\tau^2$ is required in the
first place.

Usually such an estimator is obtained from a moment-type equation \cite{ru13}.
For example, the DerSimonian--Laird \citep{dl} estimator of $\tau^2$ is
\[
\hat{\tau}^2_{\mathrm{DL}} = \frac{ \sum_i (x_i-\delta_{\mathrm{GD}})^2 s_i^{-2}
-n+1}{\sum_{i}s_i^{-2}-\sum_{i}s_i^{-4}/ (\sum_{i}s_i^{-2}
)}
\]
with $\delta_{\mathrm{GD}} =\sum_i s_i^{-2}x_i /\sum_{i} s_i^{-2}$
denoting the Graybill--Deal
estimator of $\mu$.
The popular DerSimonian--Laird $\mu$-estimator is obtained from (\ref
{es}) by using the positive part of
$ \hat{\tau}^2_{\mathrm{DL}}$.

Similarly the estimator of $\tau^2$,
\[
\hat{\tau}^2_{\mathrm{H}}=\frac{\sum_i(x_i -\bar{x})^2 -(n-1)\sum_i
s_i^2/n}{n-1},
\]
leads to the Hedges estimator of $\mu$.

The paper questions the wisdom of using under all circumstances the
tradition of plugging in $\tau^2$ estimators
to get $\mu$ estimators. Indeed the routine of plug-in estimators may
lead to poor
procedures. For example, by replacing the unknown $\tau^2$ by $\hat
{\tau}^2$ in the above formula for
$ \operatorname{Var}(\xt)$,
one can get a flagrantly biased estimator which leads to inadequate
confidence intervals for $\mu$.

A large class of weighted means statistics is motivated by the form of
Bayes procedures derived in Section~\ref{qe}.
These statistics which typically \emph{do not} admit the
representation (\ref{es})
induce estimators of the weights (\ref{we}) which shows the primary
role of $\mu$-estimation.

The main results of this work are based on a canonical representation
of the restricted likelihood function in terms
of independent normal random variables and possibly of some $\chi
^2$-random variables. An important relationship
between the weighted means statistics with weights of the form (\ref
{we}) and linear combinations of $x$'s, which are
shift invariant and independent, follows from this fact.
Our representation transforms the original problem to that of
estimating curve-confined expected values of independent
heterogeneous $\chi^2$-random variables. This reduction makes it possible
to describe the risk behavior of the weighted means statistics whose
weights are determined by a quadratic form.

We make use of the concept of permissible estimators which cannot be
uniformly improved in terms of the
differential inequality in Section~\ref{inad}.
This inequality shows that the sample mean exhibits the Stein-type
phenomenon being an inadmissible estimator
of $\mu$ under the quadratic loss when $n >3$.
A risk function for the weights in a weighted means statistic whose main
purpose is $\mu$-estimation is suggested in Section~\ref{rrisk}.
It is shown there that under this risk the sample mean is not even minimax.
Section~\ref{mini} discusses the case of approximately equal
uncertainties, and Section~\ref{exa} gives an example.
The derivation of the canonical representation of the likelihood function
is given in the \hyperref[app]{Appendix}; the proof of Theorem~\ref{th1} is delegated to the
Electronic Supplement~\cite{supp}.

\section{Estimating the common mean}

\subsection{Restricted likelihood, heterogeneity variance estimation
and quadratic forms}\label{ba}

The setting with the common mean $\mu$ and the heterogeneity variance
$\tau^2$ described in Section~\ref{in}
is a special case of a mixed linear model where statistical inference
is commonly based on
the restricted (residual) likelihood function.

The (negative) restricted log-likelihood function (\cite{scm}, Section~6.6)
has the form
\[
\mathcal{L}=\frac{1}{2} \biggl[ \sum_i
\frac{(x_i-\xt)^2}{ \tau
^2+s_i^2}+\sum_i\log\bigl(
\tau^2+s_i^2\bigr)+ \log \biggl(\sum
_i \frac{1}{\tau^2+s_i^2} \biggr) \biggr].
\]
It is possible that some of $s_i^2$ are equal; let $s_i^2$ have the
multiplicity $\nu_i, \nu_i \geq1$, so
that $\sum\nu_i=n$. Then with the index $i$ now taking values from
$1$ to $ p$,
\[
\mathcal{L}=\frac{1}{2} \biggl[ \sum_i
\frac{\nu_i(\bar{x}_i-\xt
)^2}{\tau^2+s_i^2} +\sum_i\nu_i\log
\bigl(\tau^2+s_i^2\bigr)+\log \biggl(\sum
_i \frac{\nu
_i}{\tau^2+s_i^2} \biggr) +\sum
_{i}\frac{(\nu_i-1)
u_i^2}{ \tau^2+s_i^2} \biggr].
\]
Here, $p$ denotes the number of pairwise different $s_i^2$, $\bar{x}_i
=\sum_{k: s_k=s_i} x_k/\nu_i$
represents the average of $\nu_i$
$x$'s corresponding to the particular $s_i^2, i=1, \ldots, p$, and
$u_i^2$ is their sample variance when
$\nu_i \geq2$. To simplify the notation, we write $x_i$ for $\bar{x}_i$,
so that $\xt= \sum_{i} \nu_i (\tau^2+s_i^2)^{-1} x_i/\sum_{i}\nu
_i (\tau^2+s_i^2)^{-1}$.
In our problem $x=(x_1, \ldots, x_{p})$ and $u_i^2= \sum_{k:
s_k=s_i}(x_k- \bar{x}_i)^2/(\nu_i-1)$, $\nu_i >1$,
form a sufficient statistic for $\mu$ and $\tau^2$.

Throughout this paper, we assume that $p\geq2$. Otherwise
all $\mu$-estimators reduce to the sample mean (but see Section~\ref{mini} where $\tau^2$-estimation for
equal uncertainties is considered).
The results in the \hyperref[app]{Appendix} relate the likelihood function $\mathcal
{L}$ to the joint density of $p-1$
independent normal, zero mean random variables $y_1, \ldots, y_{p-1}$.
The $(p-1)$-dimensional normal random vector $y =(y_1, \ldots, y_{p-1})^{\mathrm{T}}$
which is a linear transform of $x$ has zero mean (no matter what $\mu$
is) and the covariance matrix,
$\operatorname{diag}(\tau^2+ t_1^2, \ldots, \tau^2+ t_{p-1}^2)$,
with $t_1^2,\ldots,t_{p-1}^2$ larger than $\min s_i^2$.

To find these numbers, we introduce
the polynomial $P(v) =\prod_i (v+s_i^2)^{\nu_i}$ of degree $n$, and
its minimal annihilating polynomial
$M(v) =\prod_{i} (v+s_i^2)$ which has degree $p$. Define
%
\begin{equation}
\label{def} Q(v)=M(v) \frac{ P^\prime(v)}{P(v)}= \sum_i
\nu_i \prod_{k: k \neq
i} \bigl(v+s_k^2
\bigr).
\end{equation}
Thus $Q$ is a polynomial of degree $p-1$ which has only real (negative) roots,
denoted by $-t_1^2, \ldots,-t_{p-1}^2$ (coinciding with the roots of
$ P^\prime$ different from $-s_1^2, \ldots, -s_p^2$). Thus
$Q(v) =n \prod_{j} (v+t_j^2)$. Note that $M(-t_j^2) \neq0$. When $\nu
_i \equiv1$, $ M(v)=P(v)$, and $Q(v)=P^\prime(v)$.

According to (\ref{def}),
\[
\sum_i\log\bigl(\tau^2+s_i^2
\bigr)+\log \biggl(\sum_i \frac{\nu_i}{\tau
^2+s_i^2}
\biggr) =\sum_j \log\bigl(\tau^2+t_j^2
\bigr) +\log n,
\]
so that by using (\ref{co7}) one gets
%
\begin{eqnarray}
\label{rl} &&\mathcal{L} =\frac{1}{2} \biggl[ \sum
_j \frac{y_j^2}{\tau^2+t_j^2} +\sum_j
\log\bigl(\tau^2+t_j^2\bigr)
\nonumber
\\[-8pt]
\\[-8pt]
&&\hphantom{\mathcal{L} =\frac{1}{2} \biggl[}+\sum
_i \frac{(\nu_i-1) u_i^2}{\tau^2+s_i^2} +\sum_i(
\nu_i-1)\log \bigl(\tau^2+s_i^2
\bigr) +\log n \biggr].
\nonumber
\end{eqnarray}

The representation (\ref{rl}) of the restricted likelihood function
very explicitly takes into account one degree of freedom used for
estimating $\mu$,
as it corresponds to $p-1$ independent zero mean, normal random
variables $y_j$
with variances $\tau^2 +t_j^2, j=1, \ldots, p-1$. In addition, this\vspace*{-1pt}
likelihood includes independent $u_i^2$,
each being a multiple of a $\chi^2$-random variable with $\nu_i-1$
degrees of freedom.
When $\nu_i >1$, $u_i^2$ is an unbiased estimator of $\tau^2 +s_i^2$,
$u_i^2 \sim(\tau^2+s_i^2)\chi^2_{\nu_i-1}/(\nu_i-1)$. For $\nu_i
=1$, $u_i^2 =0$ with probability one.
According to the sufficiency principle, all statistical inference about
$\tau^2$
involving the restricted likelihood can be based exclusively on
$y_1^2,\ldots,y_{p-1}^2$ and $u_1^2,\ldots,u_p^2$.
Their joint distribution forms a curved exponential family whose
natural parameter
is formed by $(\tau^2+t_j^2)^{-1}$ (and perhaps by some $(\tau
^2+s_i^2)^{-1}$).

Evaluation of the restricted maximum likelihood estimator (REML)
$\hat{\tau}^2$ is considerably facilitated by employing $y_1^2,
\ldots, y_{p-1}^2$ and $u_1^2, \ldots, u_p^2$. Indeed
(\ref{rl}) shows that this estimator can be determined by simple
iterations as
\[
\hat{\tau}^2=\biggl({\sum_j
\frac{y_j^2-t_j^2}{(\hat{\tau
}^2+t_j^2)^2}+\sum_i\frac{(\nu_i-1)(u_i^2-s_i^{2})}{
(\hat{\tau}^2+s_i^2)^2}}
\biggr)\Big/\biggl({\sum_j\frac{1}{(\hat{\tau}^2+t_j^2)^2} +\sum
_i\frac{\nu_i-1}{(\hat{\tau}^2+s_i^2)^2}}\biggr)
\]
with $ \hat{\tau}^2_{\mathrm{DL}}$ as a good starting point, and
truncation at zero if the iteration process converges to a negative number.
Thus, $\hat{\tau}^2$ is related to a quadratic form whose coefficients
are inversely proportional to the estimated variances of $y_j^2-t_j^2$
and of $u_i^2 -s_i^2$
(cf. \citep{efm}, Section~8).\vspace*{1pt}

The form of the likelihood function $\mathcal{L}$ also
motivates the moment-type equations based on general quadratic\vspace*{-1pt} forms,
$\sum_j q_j y_j^2 +\sum_i (\nu_i-1)r_i u_i^2$ with positive
constants $q_j, r_i$.
The moment-type equation written in terms of random variables $y_1^2,
\ldots, y_{p-1}^2$ and $u_1^2, \ldots, u_p^2$ is
\[
E\biggl[ \sum_j q_j
y_j^2 +\sum_i (
\nu_i-1)r_i u_i^2\biggr] =
\biggl[ \sum_j q_j +\sum
_i (\nu_i-1)r_i\biggr]
\tau^2+ \sum_j q_j
t_j^2 +\sum_i (
\nu_i-1)r_i s_i^2.
\]
Then the estimator of $\tau^2$ by the method of moments is
\[
\hat{\tau}^2 = \frac{ \sum_jq_j (y_j^2 -t_j^2)+ \sum_i(\nu_i-1)
r_i (u_i^2- s_i^{2})}{
\sum_jq_j + \sum_i(\nu_i-1) r_i }.
\]
Unless $\tau^2$ is large, the probability that $ \hat{\tau}^2$ takes
negative values is non-negligible.
Non-negative statistics $\hat{\tau}^2_+ = \max(\hat{\tau}^2,0)$ are
used to get $\mu$-estimators of the form (\ref{es}).

The representations of two traditional statistics in Section~\ref{in}
easily follow,
\[
\hat{\tau}^2_{\mathrm{DL}}= \frac{\sum_j t_j^{-2} y_j^2+ \sum_i (\nu
_i-1)s_i^{-2} u_i^2 -n+1 }{\sum_j t_j^{-2}
+ \sum_i(\nu_i-1) s_i^{-2} }
\]
and
\[
\hat{\tau}^2_{\mathrm{H}} = \frac{ \sum_j (y_j^2 -t_j^2)+\sum_i(\nu
_i-1)(u_i^2- s_i^{2})}{n-1}.
\]

A different method-of-moments procedure suggested by Paule and Mandel
\citep{mp}
is based on solving the equation,
\[
\sum_i \frac{\nu_i (x_i-\hat{\mu})^2} { \tau^2+s_i^{2}} +\sum
_i \frac{(\nu_i-1) u_i^2}{ \tau^2+ s_i^{2} } = n-1,
\]
which has a unique positive solution, $\tau^2= \hat{\tau}^2_{\mathrm{MP}}$,
provided that
$ \sum_i \nu_i (x_i-\hat{x}_{\mathrm{GD}})^2 s_i^{-2}  +\sum_i (\nu_i-1) u_i^2 s_i^{-2} > n-1$. If this inequality does not
hold, $ \hat{\tau}^2_{\mathrm{MP}}=0$.
Because of (\ref{co7}), the equation can be rewritten in terms of
$y$'s and $t$'s as
\[
\sum_j \frac{y_j^2}{ \tau^2+t_j^{2}} +\sum
_i \frac{(\nu_i-1)
u_i^2}{ \tau^2+ s_i^{2} } = n-1.
\]
This representation allows for an explicit form of $ \hat{\tau
}^2_{\mathrm{MP}}$ in some cases.

Indeed, when $n=p=2$, $\hat{\tau}^2_{\mathrm{MP}}= \hat{\tau}^2_{\mathrm{DL}}= \hat
{\tau}^2_{\mathrm{H}}=
\max[0,y_1^2-t_1^2]$, which is also the REML. When $n=p=3$,
$y_1^2/t_1^2+ y_2^2/t_2^2\geq2$,
\[
\hat{\tau}^2_{\mathrm{MP}} =\frac{y_1^2+y_2^2}{4}-\frac{t_1^2+t_2^2}{2} +
\sqrt{ \biggl(\frac{y_1^2+y_2^2}{4} \biggr)^2+ \biggl(
\frac{t_1^2-t_2^2}{2} \biggr)^2 -\frac{(t_1^2-t_2^2)(y_1^2
-y_2^2)}{4}}.
\]

We conclude this section by noticing that
the widely used heterogeneity index $I^2$ (\citep{he}, page~117) in terms
of $y$'s and $u$'s takes the from
\[
I^2=\frac{ \sum_j t_j^{-2} y_j^2+ \sum_i (\nu_i-1)s_i^{-2} u_i^2 -n+1
}{\sum_j t_j^{-2} y_j^2+ \sum_i (\nu_i-1)s_i^{-2} u_i^2 },\qquad 0 \leq I^2 <1.
\]

\subsection{Weighted means statistics and suggested estimators}
\label{qe}

Let us look now at the generalized Bayes
estimator of $\mu$ when $\Lambda$ is a prior distribution for $\tau
^2$ while $\mu$ has the uniform
(non-informative) prior. Under the quadratic loss this estimator has
the form with
$\mathcal{L}$ given in (\ref{rl})
%
\begin{equation}
\label{rli} \delta= \frac{ \int_0^\infty \xt\exp\{-\mathcal{L}\} \,\mathrm{d}\Lambda
(\tau^2)}{ \int_0^\infty\exp\{-
\mathcal{L}\} \,\mathrm{d}\Lambda(\tau^2)} =\sum_i
\omega_i x_i.
\end{equation}
Thus $\delta$ is a weighted means statistic with
normalized weights, $\omega_i\propto\nu_i \int_0^\infty(\tau
^2+s_i^2)^{-1}[\sum\nu_k(\tau^2+s_k^2)^{-1}]^{-1}
\exp\{-\mathcal{L}\} \,\mathrm{d}\Lambda(\tau^2)$, $\sum_i \omega_i =1$,
which are shift invariant, $ \omega_i (x_1+c, \ldots, x_p+c)= \omega
_i (x_1, \ldots, x_p)$ for any real $c$.
(Any function of $y_1, \ldots, y_{p-1}$ is shift invariant.)
Indeed the use of restricted likelihood is tantamount to the practice of
weighted means statistics with invariant weights as $\mu$ estimators.
(cf. \citep{scm}, Section~9.2).

Formula (\ref{rep}) in the \hyperref[app]{Appendix} gives
\[
\delta= \bar{x}- \sum_j \frac{ \int_0^\infty(\tau
^2+t_j^2)^{-1}\exp\{-\mathcal{L}\} \,\mathrm{d}\Lambda(\tau^2)} {
\int_0^\infty\exp\{-\mathcal{L}\} \,\mathrm{d}\Lambda(\tau^2)}
\sqrt{b_j} y_j = \bar{x}- \sum
_j w_j\sqrt{b_j} y_j
\]
with $y_j$ discussed in Section~\ref{ba}.
Positive coefficients $b_j$ (the diagonal elements of the diagonal
matrix $A^{\mathrm{T}}J^{-1} A$ defined in Lemma~\ref{lem}) can
be found from (\ref{co6}) or rather from (\ref{co5}); $w_j$ is the
posterior mean of $(\tau^2+t_j^2)^{-1}$,
\[
w_j=-2\frac{\partial} {\partial y_j^2} \log\lambda\bigl(y_1^2,
\ldots, y_{p-1}^2, u_1^2,\ldots,
u_p^2\bigr)
\]
with\vspace*{-1pt} $ \lambda= \int_0^\infty\exp\{-\mathcal{L}\} \,\mathrm{d}\Lambda(\tau
^2)$. Thus positive $w_j$
is designed to estimate $(\tau^2+t_j^2)^{-1}$, $w_j \leq t_j^{-2}$,
and as a function of $y_\ell^2$, $w_j$ decreases.
The inequalities, $t_j^2 < t_\ell^2$, and $w_j > w_\ell$, are equivalent.

If $p >2$ and the support of $\Lambda$ has at least two points,
$\delta$ does not admit representation
(\ref{es}) which suggests a more general class of $\mu$-estimators.
Namely, we propose to use weighted means statistics $\delta=\sum_i
\omega_i x_i $
with weights $\omega_i = 1/p - \sum_{j} w_j A_{ij}$.
The Bayes weights belong to a smaller part of this polyhedron, namely
to the convex hull
of the vectors with coordinates $ (\tau^2+t_1^2)^{-1},\ldots,(\tau
^2+t_{p-1}^2)^{-1}$ for $\tau^2 \geq0$.
If $\hat{\tau}^2$ is an estimate of $\tau^2$, the weights
corresponding to (\ref{es}),
%
\begin{equation}
\label{weig} w_j= \frac{1}{\hat{\tau}^2+t_j^2 },
\end{equation}
lie on the boundary of this convex hull. A corner point,
$(t_1^{-2}, \ldots, t_{p-1}^{-2})$, of the convex hull always is an
inner point of the polyhedron.

Thus the focus in this paper is on estimators $\delta$ of $\mu$,
which admit the representation,
%
\begin{equation}
\label{ne} \delta=\sum_i
\omega_i x_i = \bar{x}- \sum
_j \sqrt{b_j} w_j
y_j
\end{equation}
with $w_j, 0\leq w_j \leq t_j^{-2}, y_j$ and $b_j$ as defined above.
The last term in the right-hand side of (\ref{ne}) can be viewed as an
arguably necessary heterogeneity correction
to $\bar{x}$.

Notice that (\ref{ne}) does not need an estimate of $\tau^2$ as a
prerequisite.
Since $w_j$ is an approximation to $(\tau^2+t_j^2)^{-1}$, when $n=p$,
the form of the REML $\hat{\tau}^2$ in Section~\ref{ba} suggests
such an estimator:
$[\sum w_j^2 (y_j -t_j^2)]_+/\sum w_j^2$.
If some of the multiplicities exceed one, an estimator of
$(\tau^2+s_i^2)^{-1}$ can be derived from $w_j$ by replacing $t_j^2$
by $s_i^2$.
According to (\ref{rep}), $\xt$ as well as $\bar{x}$, has the form
(\ref{ne}). In fact, all traditional statistics
(\ref{es}) admit this representation.

\subsection{Estimation of multivariate normal mean and permissible procedures}
\label{inad}

We look now at the quadratic risk behavior of $\mu$-estimators of the
form (\ref{ne}).
If $\delta=\sum_i \omega_i x_i$ is such an estimator
with positive normalized weights $\omega_i$ which are shift invariant
functions of $x_1, \ldots, x_p$,
then it is unbiased. Its variance does not depend on $\mu$ and can be
written as
%
\begin{equation}
\label{qr} \operatorname{Var}(\delta)= \operatorname{Var}(\xt) + E( \delta-
\xt)^2= \biggl[\sum_i
\frac{\nu
_i}{ \tau^2+s_i^2 } \biggr] ^{-1} + E( \delta-\xt)^2
\end{equation}
by independence of $\xt$ and $\delta-\xt$.
This and more general decompositions of the mean squared error are
discussed by Harville \citep{ha}.
The second term in the right-hand side of this identity is an important
variance component
which shows how well $\delta$ approximates the optimal but unavailable
$\xt$,
and which relates our setting to the classical estimation problem of the
multivariate $(p-1)$-dimensional normal mean.
%
\begin{prp}\label{pro2.1}
If the coefficients $w_j = w_j(y_1^2, \ldots, y_{p-1}^2, u_1^2,
\ldots, u_p^2)$
defining the estimator (\ref{ne}) are piecewise differentiable in
$y$'s, then
\begin{eqnarray*}
\operatorname{Var}(\delta) &=&\operatorname{Var}(\xt) + \sum
_j b_jE y_j^2 \biggl(
w_j-\frac{1}{\tau
^2 +t_j^2} \biggr) ^2
\\
&=&\operatorname{Var}(\bar{x})+ E\sum_j
b_j \biggl(f_j^2-2\frac{\partial}{\partial
y_j}
f_j \biggr),
\end{eqnarray*}
where $f_j =y_j w_j$. When $p >3$, $\bar{x}$ is an inadmissible
estimator of $\mu$ under the quadratic loss.
\end{prp}

The omitted proof of Proposition~\ref{pro2.1} is based on (\ref{rep}), (\ref{co5}),
and on familiar integration by parts technique. It demonstrates linkage
of our situation
to the differential inequality of a statistical estimation problem
\citep{br}. Namely, if for some vector
$\theta$, $ Y \sim N_{p-1}(\theta, I)$, then
$\sum_j b_j ( f_j^2 -2\,\partial f_j/\partial y_j )$ is an unbiased estimate
of $\sum_j b_j E (Y_j +f_j(Y) -\theta_j)^2 -\sum_j b_j \theta_j^2$.
Therefore $Y+g(Y)$, $g=(g_1, \ldots, g_{p-1})^{\mathrm{T}}$,
improves on $Y+f(Y)$, $f=(f_1, \ldots, f_{p-1})^{\mathrm{T}}$, as a $\theta$-estimator
provided that for all values $Y_1, \ldots, Y_{p-1}$,
%
\begin{equation}
\label{per} \sum_j b_j \biggl(
f_j^2 -2 \frac{\partial} {\partial Y_j} f_j \biggr)
\geq \sum_j b_j \biggl(
g_j^2 -2\frac{\partial} {\partial Y_j} g_
j \biggr).
\end{equation}
Following \citep{ru7}, let us call a (piecewise differentiable) vector
function $f
$ \emph{permissible} if (\ref{per}) does not have any solutions $g$
providing a strict inequality at some point.
Thus, $Y+f$ is a permissible estimator of the vector normal mean
$\theta$ if and only if the corresponding
scalar $\mu$-estimator, $\bar{x}- \sum_j\sqrt{b}_j f_j$, cannot be
improved upon in the sense of differential
inequality (\ref{per}). Since for $p >3$, $f\equiv0$ is not a
permissible function,
the sample mean $\bar{x}$ is inadmissible in the original setting.
Indeed the left-hand side of (\ref{per}) is negative for $f_j^{{S}} =y_j
w_j^{{S}}, w_j^{{S}}=
(p-3)/(b_j\sum y_\ell^2/b_\ell)$ proving this statement.

The differential operator in (\ref{per}) does not involve $t$'s or
$s$'s, but
in our problem only functions $f_j$ such\vspace*{-1pt} that $|f_j| \leq|y_j| t_j^{-2}$
and $f_j/y_j \geq0$ are of interest. Since $(\tau^2+t_j^2)^{-1}$ is
positive and cannot exceed $t_j^{-2}$,
according to the first equality in Proposition~\ref{pro2.1}, $w_j$ can be improved
by $\max[0,\min(w_j, t_j^{-2})]$.\vspace*{1pt}

The proof of Theorem~1 in \citep{br} shows that any permissible $w_j$
in our situation is of the form
\[
w_j=\max \biggl[0,\min \biggl(-\frac{\partial}{\partial y_j^2}\log \lambda,
\frac{1}{t_j^2} \biggr) \biggr]
\]
with some piecewise differentiable positive\vspace*{-1pt} function $\lambda =\lambda(y_1^2, \ldots, y_{p-1}^2, u_1^2, \ldots, u_p^2)$.
When $n=p$ and $\lambda=\lambda(q)$ for a positive quadratic form
$q=\sum_j q_jy_j^2, q_j >0$, one\vspace*{-1pt} gets
$w_j=\min[-q_j(\log\lambda)^\prime(q),t_j^{-2}]$.
If there are multiplicities exceeding one,
the quadratic form $ q$ is to be\vspace*{1pt} replaced by $q=\sum_j q_j y_j^2 +\sum_i(\nu_i-1)r_i u_i^2$.
For example, the function, $\lambda(q) =q^{-\alpha},\alpha>0$, leads
to the estimator
(\ref{ne}) with
%
\begin{equation}
\label{wco} w_j=\min \biggl[\frac{\alpha q_j}{\sum_\ell q_\ell y_\ell^2 +\sum_i(\nu_i-1)r_i u_i^2},
\frac{1}{t_j^2} \biggr].
\end{equation}

The statistic
$w_j^{{JS}}=\min(w_j^{{S}},t_j^{-2})$, corresponding to $ q_j =b_j^{-1},
\alpha=p-3$, when $n=p$ is similar to
the positive part of the Stein estimator of the vector normal mean
which improves over~$Y$. However, in the meta-analysis context it is
desirable having the coefficients $q_j$
of the same ordering as $t_j^{-2}$, and this condition may not hold for
$q_j \propto b_j^{-1}$.
As a matter of fact, despite\vspace*{-1pt} doing better than $w_j \equiv0$ or
$w_j^{{S}}$, the weights\vspace*{1pt} $w_j^{{JS}}$
do not produce a good estimator of $\mu$.
The same is true for many other procedures (\ref{wco}) satisfying
condition (\ref{xbar}) of
Theorem~\ref{th1} in the next section. This theorem shows that if $p \leq3
<n$, $\bar{x}$ is an inadmissible
estimator of $\mu$ although the function $f\equiv0$ is permissible then.

\subsection{$R$-risk and asymptotic optimality}
\label{rrisk}

According to (\ref{qr}) the variance of estimator (\ref{ne}) is
completely determined by the term,
$E(\delta- \xt)^2$, which can be interpreted as a cost of not knowing
$\tau^2$ when estimating $\mu$,
or as a new risk of $\delta$ viewed\vspace*{-1pt} as a procedure providing
approximations to $(\tau^2+t_j^2)^{-1}, j=1,
\ldots, p-1$. More conveniently, with $s^2=\sum_{i}\nu_i s_i^2/n$, define
\[
R\bigl(\delta,\tau^{2}\bigr)= \frac{E(\delta-\xt)^2} {\operatorname{Var}(\bar{x})- \operatorname{Var}(\xt)}= \frac{ E [ \sum_i(\omega_i - \omega_i^0) x_i ]^2} {
({\tau^2+s^2})/{n}- [\sum_i {\nu_i}/({ \tau^2+s_i^2 })
] ^{-1}}
\]
to be the $R$-risk of $\delta$.
Because of (\ref{rep}) and (\ref{ne}), the ensuing random loss
function has the form,
\[
L\bigl(\delta, \tau^2\bigr)= \frac{(\delta-\xt)^2}{
({\tau^2+s^2})/{n}-  [\sum_i {\nu_i}/({ \tau^2+s_i^2
}) ]^{-1}} =\frac{\sum_j  (w_j -{1}/({\tau^2+t_j^2}) )^2 b_j
y_j^2}{\sum_j  {b_j}/({ \tau^2+t_j^2}) }.
\]
This loss is invariant under a scale change of $y_j, \tau, t_j$ (or of
$x_i, \tau, s_i$).
For $\tau^2\to\infty$,
\[
\frac{\tau^2+s^2}{n}- \biggl[\sum_i
\frac{\nu_i}{ \tau^2+s_i^2
} \biggr]^{-1} \sim\frac{\sum_i \nu_i(s_i^2- s^2)^2}{n^2 \tau^2 },
\]
so that the normalizing factor in the definition of $L$ amplifies the
error in approximating $\xt$ when
$\tau^2$ is large.
The results of this section show that for estimators $\delta$
satisfying conditions of the following Theorem~\ref{th1},
\[
\operatorname{Var}(\delta)= \operatorname{Var}(\bar{x})+ \bigl[
\operatorname{Var}(\bar{x}) - \operatorname{Var}(\xt)\bigr] \bigl[R\bigl(\delta ,
\tau^{2}\bigr)-1\bigr] =\frac{\tau^2+s^2}{n} + \mathrm{O} \biggl(
\frac{1}{ \tau^2} \biggr)
\]
when $\tau^2 \to\infty$.
Thus, $\operatorname{Var}(\bar{x})=(\tau^2+s^2)/n$ is the dominating contribution to
the variance of $\delta$ when $\tau^2$ is large.
The $R$-risk is a useful tool for the comparison of estimators
(\ref{ne}), as
unlike the normalized quadratic risk, $E(\delta-\mu)^2/\operatorname{Var}(\bar
{x})$, it removes this linear in $\tau^2$ term.

If $\delta=\tx$ with an invariant $\hat{\tau}^2$,
then $R(\tx, \tau^2)$ can be interpreted as a conventional risk of
the estimator $\hat{\tau}^2$.
However under this risk large values of $\hat{\tau}^2$ are not
penalized very much no matter what $\tau^2$ is.
Indeed $\hat{\tau}^2$ is not designed to estimate $\tau^2$ itself,
but rather\vspace*{-1pt}
$(\hat{\tau}^2+t^2_j)^{-1}$ estimates $(\tau^2+t^2_j)^{-1}$ (cf.
\citep{MorrisN}, page 329).
When $n=p=2$, the estimator $\bar{x}$, which corresponds to
$\hat{\tau}^2=\infty$, is even admissible which of course cannot
happen for any
unbounded loss function.
This circumstance explains why an estimator $\hat{\tau}^2$ may have a
large quadratic risk, while
the associated estimator $\tx$ in (\ref{es}) has a small variance.
That phenomenon is known to happen in the case of the DerSimonian--Laird
procedure \citep{jbb}.

The estimator $ \bar{x}$ has a constant risk, $R(\bar{x},\tau^{2})
\equiv1$, which raises the question of its
$R$-minimaxity, i.e., if
$\inf_\delta\sup_{\tau^2} R( \delta,\tau^{2})=1$.
In contrast, for the Graybill--Deal estimator, $R(\delta_{\mathrm{GD}},\allowbreak   \tau
^2)=\tau^4
[\sum b_j (\tau^2 +t_j^2)^{-1} t_j^{-4}] /[\sum b_j (\tau^2 +t_j^2)^{-1}]$,
so that its $R$-risk, which vanishes\vspace*{-1pt} when $\tau^2=0$, grows
quadratically in $\tau^2$.
The next result gives a large class of estimators with bounded $R$-risk
improving on $\bar{x}$ when $n > 3$.

\begin{them} \label{th1}
Under notation of Section~\ref{ba},
let for $n>3$, $q=\sum_{j}q_j y_j^2+\sum_{i}(\nu_i-1)r_i u_i^2$
be a quadratic form with positive coefficients $q_j, r_i$. If $\delta$
has the form (\ref{ne})
such that for $q \to\infty, w_j \sim\alpha_j/ q, 0 < \alpha_j <
\infty$, then
%
\begin{eqnarray}
\label{isr} \lim_{\tau^2 \to\infty} R\bigl(\delta,\tau^{2}
\bigr)&=& 1- \frac{1}{\sum_j
b_j} \sum_j
b_j
\nonumber
\\
&&\phantom{1-}{}\times \biggl[2 \alpha_j E\frac{z_j^2} {\sum_\ell q_\ell z_\ell^2+
\sum_i r_i
\chi^2_{\nu_i-1}} -
\alpha_j^2 E\frac{ z_j^2} {(\sum_\ell q_\ell
z_\ell^2+
\sum_i r_i \chi^2_{\nu_i-1})^2} \biggr]\quad\quad\ \
\\
&\geq&\frac{2}{n-1},
\nonumber
\end{eqnarray}
where independent standard normal $z_1, \ldots, z_{p-1}$ are
independent of $\chi^2_{\nu_1-1},\ldots,
 \chi^2_{\nu_p-1}$.
Equal coefficients $q_j=r_i $ (and only they) provide the asymptotically
optimal quadratic form. If $q_j=r_i =1$, the optimal choice is $ \alpha
_j \equiv n-3$.
The sample mean $\bar{x}$ is not $R$-minimax, any estimator (\ref
{ne}) with weights (\ref{wco}) improves on it if
%
\begin{equation}
\label{xbar} 0 < \alpha\leq2(n-3) \frac{ \min[ \min_j q_j^2 t_j^4, \min_{i: \nu_i \geq2} r_i^2
s_i^4] \sum_j b_j q_j} {\max[ \max_j q_j^2 t_j^
4, \max_{i: \nu_i \geq2} r_i^2 s_i^4] \sum_j b_j q_j^2}.
\end{equation}
\end{them}

Theorem~\ref{th1} shows that the traditional weights (\ref{weig})
with $\hat{\tau}^2=q/\alpha$
are not asymptotically optimal unless the quadratic form $q$ coincides
(up to a positive factor) with
$q^\infty=\sum_j y_j^2+\sum_{i}(\nu_i-1) u_i^2=\sum_i\nu
_i(x_i-\bar{x})^2 +\sum_{i}(\nu_i-1) u_i^2$, and
$\alpha=n-3$.
Only then (\ref{isr}) is an equality. Thus, the Hedges estimator for
which $\alpha=n-1$ and $R(\delta_{\mathrm{H}},
\tau^{2})\sim2(n-3)^{-1}$, is not asymptotically optimal albeit its
performance is the best when $\tau^2$ is large.
For the Mandel--Paule estimator from Section~\ref{ba}, as well as for
the REML, (\ref{isr}) also holds
with the same quadratic form and the same $\alpha$.
The DerSimonian--Laird estimator is defined by the\vspace*{1pt} quadratic form
$q^0=\sum_j y_j^2/t_j^2+\sum_{i}(\nu_i-1) u_i^2/s_i^2 $ with $\alpha
= \sum_j t_j^{-2}+\sum_{i}(\nu_i-1) s_i^{-2}$.
Therefore, these three statistics are not optimal for large $\tau^2$ either.

The case when $n=p=2$ was studied in \citep{ru}. Then $\bar{x}$ is
admissible (so that it is
automatically minimax under $R$). Any estimator (\ref{ne}) has the
form (\ref{es}) with some
$ \hat{\tau}^2=  \hat{\tau}^2(y_1^2)$, and its $R$-risk grows
linearly in $\tau$,
\[
R\bigl(\delta,\tau^{2}\bigr)\sim \frac{\sqrt{2}}{ \sqrt{\uppi} } \int
_0^\infty\frac{y^2 \,\mathrm{d}y}{(\hat
{\tau}^2+s^2)^2} \tau.
\]
For $n =p= 3$, as $\tau^2 \to\infty$,
$ R(\delta,\tau^{2})\sim C \log\tau^2$ (see Electronic Supplement).
By analogy with the Stein phenomenon, admissibility of the sample mean
when $n=3$ is expected.

\subsection{Equal uncertainties and minimax value}
\label{mini}

When $n>3$,
the minimax value, $\inf_\delta\sup_{\tau^2} R( \delta,\tau
^{2})$, (which does not exceed one since
$R(\bar{x},\tau^{2}) \equiv1$)
cannot be smaller than $2 (n-1)^{-1}$. Indeed for any estimator $\delta$,
\[
\sup_{\tau^2} R\bigl( \delta,\tau^{2}\bigr) \geq
\limsup_{\tau^2 \to\infty
} R\bigl( \delta,\tau^{2}\bigr) \geq
\frac{2}{n-1}.
\]
This fact can be proven by constructing a sequence of proper prior
densities for $\tau^2$
such that the corresponding sequence of the Bayes $R$-risks converges
to $2 (n-1)^{-1}$.

Thus for large $\tau^2$, the estimators (\ref{ne}) with $q^\infty$,
$\alpha=n-3$, cannot be improved upon.
The most natural of these statistics, say, $\delta_1$ has the form
(\ref{ne}) with
%
\begin{equation}
\label{del} w_j^1 = \min \biggl(\frac{n-3}{q^\infty},
\frac{1}{t_j^2} \biggr).
\end{equation}
Another\vspace*{-1pt} modified Hedges estimator, $\delta_{\mathrm{mH}}$, has the form (\ref
{es}) with $\hat{\tau}^2= (n-3)^{-1}[q^\infty-\sum_j t^2_j-\sum_i (\nu_i-1)s_i^2]_+$
and also is asymptotically optimal although in general its performance
is worse than that of (\ref{del}).

If $\sum_i \nu_i(s_i^2- s^2)^2 \to0$, so that all $t_j^2$ and
$s_i^2$ tend to $s^2
=\sum_i \nu_i s_i^{2}/n$, $w_j(v) \to w(v)$,
\[
R\bigl( \delta,\tau^2\bigr)\to \bigl(\tau^2+s^2
\bigr)^2 \int_0^\infty \biggl[ w(v)-
\frac{1}{\tau^2+s^2} \biggr]^2 \,\mathrm{d}G_{n+1} \biggl(
\frac{v} {\tau^2+s^2} \biggr).
\]
Here and further $G_{k}$ is the distribution function of $\chi^2$-law
with $k$ degrees of freedom.
Thus if $s_i^2\approx s^2$, our problem is that of estimation of the
reciprocal of the scale parameter
$\sigma=\tau^2+s^2$ under the restriction, $\sigma\geq s^2$.
The ``data'' $v$ in this problem is $\chi^2$-distributed, $v \sim
\sigma\chi^2_{n+1}$,
and the invariant loss function, $\sigma^2( w- \sigma^{-1})^2$,
corresponds to the $R$-risk.
Then the minimax value, $2 (n-1)^{-1}$, is the same as in the
non-restricted ($s=0$) parameter case \citep{MS0}.
As in unrestricted scale parameter estimation, the generalized prior,
$\mathrm{d}\sigma/\sigma$, $\sigma\geq s^2$,
or $\mathrm{d}\tau^2/(\tau^2+s^2)$, provides a least favorable distribution.
See also \citep{MS1} for more general results.

Thus in meta-analysis problems with $s_i^2$ exhibiting little variation,
the minimax value is expected to stay close to $2 (n-1)^{-1}$. Indeed
when $w(v)=  \min(\alpha v^{-1},  s^{-2})$, $\xi=\alpha s^2/(\tau^2+s^2)$,
%
\begin{equation}
\label{rar} R\bigl( \delta,\tau^2\bigr)\to 1- \biggl(1-
\frac{\tau^4}{s^4} \biggr) G_{n+1}(\xi) - \frac{2\alpha(1-G_{n-1}(\xi))}{n-1} +
\frac{\alpha
^2(1-G_{n-3}(\xi))}{(n-1)(n-3)}.
\end{equation}
The formula (\ref{rar}) shows that the estimator (\ref{del}) is
minimax unlike $\delta_{\mathrm{mH}}$ for which
$w(v)=\{ [v- (n-1)s^2]_+ /(n-3) + s^2 \}^{-1}$.

The DerSimonian--Laird rule, $w_{\mathrm{DL}}(v)=\{ [v- (n-1)s^2]_+ /(n-1) +s^2 \}
^{-1}, \alpha=n-1$,
coincides in this situation with the REML and the Hedges estimator. For
the proper maximum likelihood estimator of $(\tau^2+s^2)^{-1}$,
$\alpha=n+1$. None of these procedures is
minimax which indicates that their good properties in meta-analysis may
be attributable to
a large number of individual studies (large $n$) or to lack of interest
in high heterogeneity (small $\tau^2$).

\begin{figure}

\includegraphics{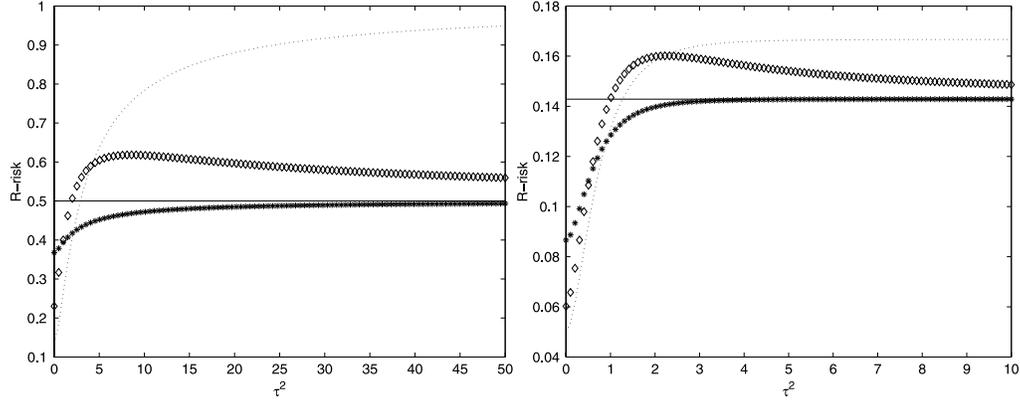}

\caption{Plots of $R$-risks of estimators corresponding to $\delta
_{\mathrm{DL}}$ (dash-dotted line),
$\delta_{\mathrm{mH}}$ (line marked by diamonds) and $\delta_1$ (line marked
by $*$) when $n=5$ (left panel),
and of the same risks when $n=15$ (right panel). The straight line
depicts the minimax value $2(n-1)^{-1}$.}\label{fi1}
\end{figure}

Figure~\ref{fi1} shows the graphs of the $R$-risk in (\ref{rar})
when $s^2=1$.
It suggests that the estimator $\delta_1$ performs quite well for
small/medium $n$'s.
Indeed this estimator is better than other procedures except for small
$\tau^2$
in which case $\delta_{\mathrm{DL}}$ dominates $\delta_{\mathrm{mH}}$ (at the price of
higher risk for larger
values of $\tau^2$).

\section{Example: $p=2$}
\label{exa}

When there are only two different values $s^2_1$ and $s_2^2$
with multiplicities $\nu_1$ and $\nu_2$, $n=\nu_1+\nu_2 >3$,
\begin{eqnarray*}
t_1^2&=&t^2= \frac{ \nu_2 s_1^2 +\nu_1 s_2^2}{n},
\\
A_{11}&=&-A_{21}= \frac{ \nu_1 \nu_2(s_1^2- s_2^2)}{n^2},
\\
b_1&=&b= \frac{ \nu_1 \nu_2(s_1^2- s_2^2)^2}{n^3},
\end{eqnarray*}
and if $s_1^2 > s_2^2$,\vspace*{-2pt}
\[
y_1=y=\frac{\sqrt{\nu_1 \nu_2} (x_1-x_2)}{\sqrt{n}}.
\]
Any estimator (\ref{ne}) has the form (\ref{es}) for some $\hat{\tau}^2$,\vspace*{-2pt}
\[
\delta=\frac{x_1+x_2}{2}-\frac{\nu_1 \nu
_2(s_1^2-s_2^2)(x_1-x_2)}{n^2 (\hat{\tau}^2+t^2)}.
\]

For $\delta_1$, $\hat{\tau}^2=[(n-3)^{-1}q^\infty-t^2]_+, q^\infty
=y^2+(\nu_1-1)u_1^2+(\nu_2-1)u_2^2$.
The modified Hedges estimator $\delta_{\mathrm{mH}}$
with $ \hat{\tau}^2=(n-3)^{-1}[q^\infty-t^2-(\nu_1-1)s_1^2-(\nu
_2-1)s_2^2 ]_+$ typically
has its $R$-risk at $\tau^2=0$ larger than that of $\delta_1$. (The
exact condition for $\delta_{\mathrm{mH}}$ to have a smaller
$R$-risk at $\tau^2=0$ than $\delta_{1}$ is: $n \geq5$, and if $\nu
_1 \leq\nu_2$, then
$\nu_1 \leq n(n-4)/(2n-5)$, $[n(n-1)-\nu_1(2n-5)] s_2^2 \geq
[n(n-4)-\nu_1(2n-5)] s_1^2$.)

The $R$-risk of $\delta_1$ at $\tau^2=0$ can be larger than
$2/(n-1)$. Indeed\vspace*{-1pt}
\[
R(\delta_1, 0) =\int_{ (n-3)}^{\infty} \biggl(
\frac{n-3}{v} -1 \biggr)^2 \,\mathrm{d}F\bigl(t^2 v
\bigr),
\]
where $F(v)$ is the distribution function of $t^2 \chi^2_3+ s_1^2 \chi
^2_{\nu_1-1} +s_2^2 \chi^2_{\nu_2-1}$.
With $ a=  t^2/[t^6 s_1^{2(\nu_1-1)}\*s_2^{2(\nu_2-1)}]^{1/(n+1)}$,
according to the Okamoto inequality \citep{MO},
$F(t^2 v) \leq G_{n+1}(a v)$. Thus,
since $[(n-3)v^{-1}-1]^2$ is an increasing function of $v, v \geq n-3$,\vspace*{-1pt}
\begin{eqnarray*}
R(\delta_1, 0) &>& \int_{ (n-3)}^{\infty}
\biggl( \frac{n-3}{v} -1 \biggr)^2 \,\mathrm{d}G_{n+1}(a
v)
\\
&=& 1-G_{n+1}\bigl(a(n-3)\bigr)
\\
&&{}-\frac{2(n-3)a[1-G_{n-1}(a(n-3))]}{n-1}
\\
&&{}+ \frac
{(n-3)a^2[1-G_{n-3}(a(n-3))]}{n-1}.
\end{eqnarray*}
This inequality shows that $ R(\delta_1, 0) \geq2(n-1)^{-1}$,
if $a <a_0<1$, where $a_0=a_0(n)$ is monotonically increasing to $1$ in
$n,a_0(4)=0.637\ldots, a_0(10)=0.798\ldots$
For small $a$, $\delta_1$ cannot have its risk at the origin smaller
than $2(n-1)^{-1}$.
For example, when $n=4, \nu_1=1, \nu_2=3$,
$R(\delta_1, 0)\leq2(n-1)^{-1}$ if and only if $s_1^2/s_2^2 \geq
0.173\ldots$\,, i.e., iff $a \geq0.679\ldots$\,.

The DerSimonian--Laird estimator $\delta_{\mathrm{DL}}$ with $\hat{\tau
}^2_{\mathrm{DL}}=(q^0-n+1)_+/(1/t^2+(\nu_1-1)/s_1^2+
(\nu_2-1)/s_2^2), q^0=y^2/t^2+(\nu_1-1)u_1^2/s_1^2+(\nu_2-1)u_2^2/s_2^2$,
has its $R$-risk at $\tau^2=0$ of the form
\[
R(\delta_{\mathrm{DL}},0)=\int_{(n-1)}^{\infty} \biggl[
\frac{1}{1+\kappa
^{-1}(v-n+1)}-1 \biggr]^2 \,\mathrm{d}G_{n+1}(v)
\]
with $\kappa= 1+ t^2[(\nu_1-1)/ s_1^2 +(\nu_2-1)/s_2^2]$.

For the estimator $\delta_0$ defined by (\ref{ne}),
\[
w_j^0 = \min \biggl(\frac{n-1}{q^0}, 1 \biggr)
\frac{1}{t_j^2},
\]
so that $\hat{\tau}^2=t^2[q^0/(n-1)-1]_+$. Its risk at $\tau^2=0$,
\[
R(\delta_{0}, 0)= \int_{(n-1)}^{\infty}
\biggl(\frac
{n-1}{v}-1 \biggr)^2 \,\mathrm{d}G_{n+1}(v),
\]
is always smaller than that of $\delta_1$.

But $\delta_0$ is also competitive against $\delta_{\mathrm{DL}}$. Indeed $
R(\delta_{0},0)< R (\delta_{\mathrm{DL}}, 0)$
if and only if $\kappa < n-1$, that is, iff
\[
\frac{(\nu_2-1)s_1^2}{\nu_1 s_2^2}+\frac{(\nu_1-1)s_2^2}{\nu_2
s_1^2} < \frac{\nu_1 \nu_2(n- 2)}{n}-2 +
\frac{n}{ \nu_1 \nu_2}.
\]
Thus provided that $\nu_1, \nu_2 >1, s_1^2/s_2^2 \approx\sqrt{ (\nu
_1-1) \nu_1/[ (\nu_2-1) \nu_2]}$,
$\delta_{0}$ improves upon $\delta_{\mathrm{DL}}$ for small $\tau^2$.
If, say, $\nu_1=1$, this domination means that $s_1^2 < s_2^2$.
Thus, when one study reports a smaller uncertainty than all other
studies whose standard errors are approximately
equal, $\delta_0$ improves upon the DerSimonian--Laird estimator for
small $\tau^2$.

However, there is no uniform domination as the condition, $\kappa <
n-1$, means that for large $\tau^2,
R(\delta_0,\tau^{2})> R(\delta_ {\mathrm{DL}}, \tau^{2})$.

\section{Conclusions}

Author's attempt was to give a perspective of a meta-analysis setting
from the point of view of the statistical
decision theory. Although concepts like admissibility or minimaxity
have not so far generated much interest
among meta-analysts, there is a realization that different desirable
qualities of the employed procedures call for
different loss functions. The quadratic loss for the mean effect
estimators from a wide class leads in a natural way
to the $R$-risk suggested and studied in this paper. This risk strongly
recommends against the use of the sample mean
as a consensus estimate which still happens in some collaborative studies.

Moreover, the $R$-risk questions well recognized
excellent properties of the DerSimonian--Laird estimator $\delta_{\mathrm{DL}}$
in the situation when $s_i$ are almost equal,
or when one study claims a high precision while all other studies report
larger uncertainties which are about the same.
The unsatisfactory performance of the Graybill--Deal estimator $\delta
_{\mathrm{GD}}$ is well known in the latter case.
It is of interest that $\delta_0$ improves on the DerSimonian--Laird
estimator for moderate/small
$\tau^2$.
Inference on the overall effect can be obtained before the
heterogeneity variance is estimated,
but even in the simplest cases considered here there is no unique rule
dominating all others.

This paper is dedicated to the memory of George Casella
who was always interested in implications of the statistical decision
theory results
to practical estimation problems \citep{maac}.

\begin{appendix}\label{app}
\section*{Appendix}

\subsection{Partial fraction decomposition and weighted means}

Let $e$ denote unit coordinates vector whose dimension is clear from
the context,
and put $J=\operatorname{diag}(\nu_1, \ldots, \nu_p)$,
$S = \operatorname{diag}(s_1^2/\nu_1, \ldots, s_p^2/\nu_p)$.
In the used here notation of Section~\ref{ba}
the vector $x$ has the diagonal covariance matrix,
$C =\tau^2 J^{-1}+S$.

\begin{lem}\label{lem}
For any $v$ different from $-t_j^2, j=1, \ldots, p-1$, and for any
$i=1, \ldots, p$,
%
\begin{equation}
\label{coo} \frac{\nu_i }{v+s_i^2} \biggl[\sum_k
\frac{\nu_k }{v+s_k^2 } \biggr] ^{-1}= \frac{\nu_i }{n} -\sum
_{j} \frac{ A_{ij}}{ v+t_j^2},
\end{equation}
where
%
\begin{equation}
\label{coe} A_{ij}= \frac{ \nu_i M(-t_j^2)}{ Q^{\prime}(-t_j^2)(t_j^2-s_i^2)}.
\end{equation}
For any $j, j=1, \ldots, p-1$,
%
\begin{equation}
\label{co6} b_j=\sum_i
\frac{A^2_{ij}}{\nu_i}= \frac{1}{ t_j^2} \sum_i
\frac{
s_i^2 A^2_{ij}}{\nu_i}=- \frac{ M(-t_j^2)}{ Q^{\prime}(-t_j^2)}.
\end{equation}
If the $p \times(p-1)$ matrix $A$ is determined by
its elements $A_{ij}$ in (\ref{coe}), then
%
\begin{equation}
\label{co0} A^{\mathrm{T}}  e =0,
\end{equation}
and
%
\begin{equation}
\label{co2} A  e =\frac{1}{n} J\bigl(S-s^2 I\bigr)
 e ,\qquad s^2=\frac{\sum_i \nu_i s_i^{2}}{n}.
\end{equation}
The matrices $ A^{\mathrm{T}} J^{-1} A = \operatorname{diag}(b_1, \ldots, b_{p-1})$
and $ A^{\mathrm{T}} S A = \operatorname{diag}(b_1t_1^2, \ldots, b_{p-1}t_{p-1}^2)$ are
diagonal, and
%
\begin{equation}
\label{co3} A\bigl(A^{\mathrm{T}} J^{-1} A\bigr)^{-1}A^{\mathrm{T}}
= J -\frac{J e   e ^{\mathrm{T}} J}{ e ^{\mathrm{T}} J  e }.
\end{equation}
With $\h= ((\tau^2 +t_1^2)^{-1}, \ldots,(\tau^2 +t_{p-1}^2)^{-1})^{\mathrm{T}}$,
%
\begin{equation}
\label{co10} A^{\mathrm{T}} J^{-1}C^{-1} J^{-1}A
= \operatorname{diag}\bigl( A^{\mathrm{T}} J^{-1}A \h\bigr) + \biggl(
\sum_i \frac{\nu_i} { \tau^2 +s_i^2} \biggr)
\bigl(A^{\mathrm{T}} J^{-1}A \h\bigr) \bigl(A^{\mathrm{T}}
J^{-1}A \h\bigr)^{\mathrm{T}},
\end{equation}
and
%
\begin{equation}
\label{co7} \sum_i \frac{\nu_i (x_i-\xt)^2}{ \tau^2+s_i^2}= \sum
_j \frac
{y_j^2}{\tau^2+t_j^2}.
\end{equation}
\end{lem}

\begin{pf}
By the definition of the polynomial $Q$ in Section~\ref{ba},
\begin{eqnarray*}
\frac{\nu_i }{n} -\frac{\nu_i }{v+s_i^2} \biggl[\sum_k
\frac{\nu
_k }{v+s_k^2 } \biggr] ^{-1} &=& \frac{\nu_i }{n} -
\frac{\nu_i P(v) }{(v+s_i^2) P^\prime(v)}
\\
&=& \frac{\nu_i }{n} -\frac{\nu_i M(v) }{(v+s_i^2) Q (v)}= \frac{ \nu_i [ \prod_j (v+t_j^2) - \prod_{k \neq i} (v+s_k^2) ]}{Q(v)}
\end{eqnarray*}
with the right-hand side of this identity being the ratio of two polynomials
of degree $p-2$ and $p-1$, respectively.
The formulas (\ref{coo}) and (\ref{coe}) easily follow from the
classical result on
partial fraction decomposition for such ratios.

For any fixed $j$,
\[
\sum_i \frac{\nu_i}{s_i^2-t_j^2} =
\frac{ P^\prime(-t_j^2)}{P(-t_j^2)}=0,
\]
so that (\ref{co0}) follows from (\ref{coe}),
\[
\sum_i A_{ij}= \frac{ M(-t_j^2)}{ Q^{\prime}(-t_j^2)} \sum
_i \frac
{\nu_i}{s_i^2-t_j^2} =0.
\]

By equating coefficients at $v^{n+p-2}$ of $QP$ and $M P^\prime$, one\vspace*{-1pt}
gets $\sum t_j^2=\sum s_i^2-s^2$.
The comparison of coefficients at $v^{p-2}$ of two equal polynomials,
$n\sum_j A_{ij} \prod_{\ell\neq j} (v+t_\ell^2) $ and $\nu_i
[\prod_j (v+t_j^2) -
\prod_{k \neq i} (v+s_k^2)]$, shows that
\[
n \sum_j A_{ij} = \nu_i
\biggl( \sum_j t_j^2-
\sum_{k \neq i} s_k^2 \biggr) =
\nu_i \biggl(s_i^2 - \sum
_{k} \nu_k s_k^2/n
\biggr),
\]
which implies (\ref{co2}).

For any $j$
\[
\sum_i \frac{\nu_i}{(s_i^2-t_j^2)^2} =-
\frac{Q^{\prime}
(-t_j^2)}{M(-t_j^2)},
\]
and
\[
\sum_i \frac{\nu_is_i^2}{(s_i^2-t_j^2)^2} =
t_j^2 \sum_i
\frac{\nu
_i}{(s_i^2-t_j^2)^2} = - t_j^2 \frac{ Q^{\prime} (-t_j^2)}{M(-t_j^2)},
\]
so that (\ref{co6}) is established by substituting (\ref{coe}) for $A_{ij}$.

For any different $j$ and $\ell$
\[
0=\sum_i \frac{\nu_i}{ s_i^2-t_j^2} -\sum
_i \frac{\nu
_i}{s_i^2-t_\ell^2} = \bigl(t_j^2-t_\ell^2
\bigr) \sum_i \frac{\nu_i}{ (s_i^2-t_j^2) (s_i^2-t_\ell
^2)},
\]
which implies that $ \sum_i\nu_i^{-1} A_{ij} A_{i\ell}=0$, or that $
A^{\mathrm{T}} J^{-1} A $ is a diagonal matrix.

This argument also shows that
\[
\sum_i \frac{s_i^2 A_{ij} A_{i\ell}}{\nu_i}=0,
\]
as
\[
\sum_i \frac{\nu_is_i^2}{ (s_i^2-t_j^2) (s_i^2-t_\ell^2)} =t_j^2
\sum_i \frac{\nu_i}{ (s_i^2-t_j^2) (s_i^2-t_\ell^2)}=0.
\]

To prove (\ref{co3}),
observe that for $i \neq k$, the $(i,k)$th element of the matrix
$A(A^{\mathrm{T}} J^{-1} A)^{-1}A^{\mathrm{T}} $ has the form,
\begin{eqnarray*}
\sum_ {j} \frac{A_{i j}A_{k j}}{b_j} &=&
\nu_i \nu_k \sum_ {j}
\frac
{b_j}{ ( t_j^2- s_i^2) (t_j^2-s_k^2)}
\\
&=& \frac{\nu_i \nu_k}{ s_i^2- s_k^2} \biggl[ \sum_ {j}
\frac{b_j}{
t_j^2-s_i^2} - \sum_ {j} \frac{b_j}{ t_j^2-s_k^2}
\biggr]=- \frac{\nu_i \nu
_k}{n}.
\end{eqnarray*}
Here we used the facts that $A_{ij}= -\nu_i b_j/ ( t_j^2- s_i^2)$, and $
\sum_{ j}b_{j}(t_j^2-s_i^2)^{-1} = ( s^2-s_i^2)/n$.

To determine the diagonal elements of $A(A^{\mathrm{T}} J^{-1} A)^{-1}A^{\mathrm{T}} $,
observe that
according to the definition of $Q,Q(-s_i^2)/ M^{\prime} (-s_i^2)=\nu_i$.
Therefore for any $i$,
\[
\sum_ {j} \frac{A_{i j}^2}{b_j} = -\nu_i
\sum_ {j} \frac{A_{i j}}{
t_j^2- s_i^2} =
\nu_i^2 \biggl[ \frac{M^{\prime} (-s_i^2)}{Q(-s_i^2)}-\frac{1}{n}
\biggr] = \nu_i -\frac{ \nu_i^2 }{n}.
\]
Thus, (\ref{co3}) holds.

Because of (\ref{coe}),
%
\begin{equation}
\label{co9} A^{\mathrm{T}} J^{-1}C^{-1}  e  = -
\biggl(\sum_i \frac{\nu_i} { \tau^2 +s_i^2} \biggr)
\bigl(A^{\mathrm{T}} J^{-1}A\bigr) \h.
\end{equation}

To prove (\ref{co10}) for fixed $j,\ell$, multiply (\ref{coo}) by
$A_{i j}, A_{i \ell}$, divide by $\nu_i^2$,
and sum up over $i$ to get the following expression for the $(j, \ell
)$th element of the matrix
$ A^{\mathrm{T}} J^{-1}C^{-1} J^{-1}A $,
\[
\biggl(\sum_i \frac{\nu_i} { \tau^2 +s_i^2} \biggr) \biggl[
\frac{ \delta_{j \ell} b_j}{n} - \sum_ {i, m} \frac{A_{i
j}A_{i \ell} A_{im}}{\nu_i^2
(\tau^2 +t_m^2)}
\biggr],
\]
where $\delta_{ij}$ is the Kronecker symbol ($\delta_{ij}=1$, if
$i=j; =0$ otherwise).
It is easy to see that $\sum_ {i} A_{i j} A_{i \ell} A_{im}\nu
_i^{-2}=0$, unless there are at least two equal
indices among $j,\ell,m$. When all three of these indices coincide,
\begin{eqnarray*}
\sum_ {i} \frac{A_{i j}^3}{\nu_i^2} &=& -
\frac{ M^3(-t_j^2)} {
[Q^{\prime} (-t_j^2)]^3} \sum_i \frac{\nu_i}{ (s_i^2-t_j^2)^3 }
\\
&=&-\frac{
M(-t_j^2)[ Q^{\prime\prime} (-t_j^2) M(-t_j^2)- 2 Q^\prime (-t_j^2)
M^\prime (-t_j^2)]
}{2 [Q^\prime (-t_j^2)]^3}
\\
&=&\frac{ b_j [ Q^{\prime\prime} (-t_j^2) M(-t_j^2)- 2 Q^\prime (-t_j^2)
M^\prime (-t_j^2)] }{2 [Q^\prime (-t_j^2)]^2} =- \frac{ b_j^2 Q^{\prime\prime} (-t_j^2)+ 2 b_j M^\prime
(-t_j^2)}{2 Q^\prime (-t_j^2)}.
\end{eqnarray*}
If, say, $m=j \neq\ell$,
\begin{eqnarray*}
\sum_ {i} \frac{A_{i j}^2 A_{i \ell}}{\nu_i^2} &=& -
\frac{
M^2(-t_j^2) M(-t_\ell^2)} {
[Q^{\prime} (-t_j^2)]^2 Q^{\prime}(-t_\ell^2)} \sum_i \frac{\nu
_i}{ (s_i^2-t_j^2)^2 (s_i^2-t_\ell^2)}
\\
&=&- \frac{ M^2(-t_j^2) M(-t_\ell^2)}{ [Q^{\prime} (-t_j^2)]^2
Q^{\prime}(-t_\ell^2)
(t_j^2-t_\ell^2)} \sum_i
\frac{\nu_i}{ (s_i^2-t_j^2)^2 }
\\
&=& \frac{ M(-t_j^2) M(-t_\ell^2)} { Q^{\prime} (-t_j^2) Q^{\prime
}(-t_\ell^2)
(t_j^2-t_\ell^2)} = \frac{ b_j b_\ell}{t_j^2-t_\ell^2}.
\end{eqnarray*}
The last formula shows that off-diagonal elements of $ A^{\mathrm{T}} J^{-1}C^{-1}
J^{-1}A $,
for $j \neq\ell$ have the form
\begin{eqnarray*}
&&- \biggl(\sum_i \frac{\nu_i}{\tau^2 +s_i^2} \biggr)
\biggl[ \frac{ b_j b_\ell} {(t_j^2-t_\ell^2) (\tau^2 +t_j^2)} +\frac{ b_j b_\ell} {(t_\ell^2-t_j^2)(\tau^2 +t_\ell^2)} \biggr]
\\
&&\quad= \biggl(\sum_i \frac{\nu_i}{\tau^2+s_i^2}
\biggr)\frac{b_j b_\ell
} {(\tau^2 +t_j^2) (\tau^2 +t_\ell^2)},
\end{eqnarray*}
that is, (\ref{co10}) holds for the off-diagonal elements.

We demonstrate now the equality of the diagonal elements
of matrices in (\ref{co10}). These elements for the matrix $ A^{\mathrm{T}}
J^{-1}C^{-1} J^{-1}A $ are
\[
b_j \biggl(\sum_i
\frac{\nu_i} { \tau^2 +s_i^2} \biggr) \biggl[ \frac{1} {n}-\sum
_{\ell\neq j} \frac{ b_\ell}{ (t_j^2-
t_\ell^2)(\tau^2 +t_\ell^2)} + \frac{ b_j Q^{\prime\prime}(-t_j^2)+2M^\prime (-t_j^2)}{2
Q^{\prime} (-t_j^2)(\tau^2+t_j^2)} \biggr].
\]

Define the polynomial $ Q_j$ by the formula,
\[
\frac{Q_j(\tau^2)}{Q (\tau^2)}= \sum_{\ell\neq j} \frac{ b_\ell
}{ (t_\ell^2-t_j^2)
(\tau^2 +t_\ell^2)}+
\frac{b_j Q^{\prime\prime}(-t_j^2)+2M^\prime
(-t_j^2)}{2Q^{\prime}(-t_j^2)(\tau^2 +t_j^2)}.
\]
Then the degree of $Q_j$ is $p-2$, and this polynomial is determined by
its values at $-t^2_\ell, \ell=1, \ldots, p-1$:
$ Q_j(-t_\ell^2)=b_\ell Q^{\prime} (-t_\ell^2)/(t_\ell^2-t_j^2)=
- M(-t_\ell^2)/(t_\ell^2 -t_j^2),\ell\neq j $, and $Q_j(-t_j^2)=b_j
Q^{\prime\prime}(-t_j^2)/2 + M^\prime(-t_j^2)$.
It follows that
\[
Q_j\bigl(\tau^2\bigr)= \frac{M(\tau^2)}{\tau^2+t_j^2}+
\frac{b_jQ(\tau
^2)}{(\tau^2+t_j^2)^2}-\frac{Q(\tau^2)}{n}.
\]
Indeed, the polynomial in the right-hand side has degree $p-2$. Since
\[
\lim_{ \tau^2\to-t_j^2}\frac{M(\tau^2) (\tau^2 +t_j^2) +b_j Q
(\tau^2)} {
(\tau^2 +t_j^2)^{2}}= M^\prime
\bigl(-t_j^2\bigr)+ \frac{b_j Q^{\prime\prime
} (-t_j^2)}{2},
\]
it assumes the same values as $Q_j$ at $-t^2_\ell, \ell=1, \ldots,
p-1$, which establishes (\ref{co10}).

Because of (\ref{co3}) and (\ref{coe}),
$x-\xt  e = x-\bar{x}  e + (\bar{x}-\xt)  e =
[I-  e   e ^{\mathrm{T}} J/( e ^{\mathrm{T}} J  e )]x +(\h
^{\mathrm{T}} A^{\mathrm{T}} x)  e $.
Thus the quadratic form in the left-hand side of (\ref{co7}) can be
written as
\begin{eqnarray*}
&&\bigl[J^{-1}A \bigl(A^{\mathrm{T}}J^{-1} A
\bigr)^{-1} A^{\mathrm{T}} x+  e  \h^{\mathrm{T}}
A^{\mathrm{T}} x\bigr]^{\mathrm{T}} C^{-1} \bigl[J^{-1}A
\bigl(A^{\mathrm{T}} J^{-1}A\bigr)^{-1} A^{\mathrm{T}} x+
 e  \h^{\mathrm{T}} A^{\mathrm{T}} x\bigr]
\\
&&\quad=y^{\mathrm{T}} \bigl[ \bigl(A^{\mathrm{T}} J^{-1}A
\bigr)^{-1/2} A^{\mathrm{T}} J^{-1}+\bigl(A^{\mathrm{T}}
J^{-1}A\bigr)^{1/2} \h  e ^{\mathrm{T}}\bigr]
\\
&&\qquad{}\times C^{-1}\bigl[ J^{-1}A \bigl(A^{\mathrm{T}} J^{-1}A
\bigr)^{-1/2} +  e  \h^{\mathrm{T}} \bigl(A^{\mathrm{T}}
J^{-1}A\bigr)^{1/2} \bigr]y
\\
&&\quad= y^{\mathrm{T}} \operatorname{diag}(\h)y,
\end{eqnarray*}
where the second equality
follows from (\ref{co10}) and (\ref{co9}).\vadjust{\goodbreak}
\end{pf}

The following important representation for $\xt$
%
\begin{equation}
\label{rep} \xt= \bar{x}- \sum_{ i,j}
\frac{A_{i j} x_i}{ \tau^2 +t_j^2} = \bar{x}- \sum_j
\frac{\sqrt{b_j} y_j}{\tau^2 +t_j^2}
\end{equation}
is a consequence of Lemma~\ref{lem}.
Here $y_j=\sum_i A_{ij}x_i/\sqrt{b_j}$ are independent normal, zero
mean random variables with the variances
$\tau^2 +t_j^2$. Indeed the normal random vector $y= (A^{\mathrm{T}} J^{-1}
A)^{-1/2} A^{\mathrm{T}} x$
has the covariance matrix $(A^{\mathrm{T}} J^{-1} A)^{-1/2}
A^{\mathrm{T}}\! C A (A^{\mathrm{T}} J^{-1} A)^{-1/2}
= \operatorname{diag}( \tau^2 +t_1^2, \ldots, \tau^2 +t_{p-1}^2)$.
Since $Ey_j(\xt-\mu)=0$, $\xt$ and $y_j$ are independent implying
independence of
$\xt$ and $\delta-\xt$ in Section~\ref{inad}.

The coefficients $A_{ij}$ provide a simple expression for $\operatorname{Var}(\bar
{x})- \operatorname{Var}(\xt)$. Indeed,
by dividing (\ref{coo}) by $\nu_i$ and multiplying it by $A_{i\ell
}$, one gets after summing up over all
$i$ and $\ell$ and using (\ref{co0}), (\ref{co2}),
\begin{eqnarray*}
\frac{P(\tau^2)}{ P^\prime(\tau^2)} \sum_{i, \ell} \frac
{A_{i\ell}}{ \tau^2 +s_i^2}&=&
\frac{M(\tau^2)}{ n Q(\tau^2)} \sum_{i} \frac{\nu
_{i}(s_i^2-s^2)}{ \tau^2 +s_i^2}
\\
&=& \biggl[\sum_i \frac{\nu_i}{ \tau^2+s_i^2 } \biggr]
^{-1} -\frac
{\tau^2+s^2}{n} = -\sum_{i, j}
\frac{A_{ij}^2}{\nu_i( \tau^2 +t_j^2)}.
\end{eqnarray*}
This formula can be written in the form,
%
\begin{equation}
\label{co5} \sum_{ j} \frac{ b_{j}}{ \tau^2 +t_j^2} =
\operatorname{Var}(\bar{x})- \operatorname{Var}(\xt) =\frac{\sum_{i} \nu_{i}(s^2-s_i^2)\prod_{k \neq i} (\tau^2
+s_k^2)}{n Q(\tau^2)},
\end{equation}
which provides the representation of the left-hand side of (\ref{co5})
as a ratio of two polynomials of
degree $p-2$ and $p-1$, respectively and which allows numerical
evaluation of $b$'s without
calculating $A_{ij}$.
\end{appendix}


\begin{supplement}
\stitle{Restricted likelihood representation and decision-theoretic
aspects of meta-analysis: Electronic supplement}
\slink[doi]{10.3150/13-BEJ547SUPP} 
\sdatatype{.pdf}
\sfilename{BEJ547\_supp.pdf}
\sdescription{The supplement contains the proof of Theorem \ref{th1}.}
\end{supplement}


\printhistory

\end{document}